# A Distributed GNE Seeking Algorithm Using the Douglas-Rachford Splitting Method

Yuanhanqing Huang and Jianghai Hu

*Abstract*—**We consider a generalized Nash equilibrium problem (GNEP) for a network of players. Each player tries to minimize a local objective function subject to some resource constraints where both the objective functions and the resource constraints depend on other players' decisions. By conducting equivalent transformations on the local optimization problems and introducing network Lagrangian, we recast the GNEP into an operator zero-finding problem. An algorithm is proposed based on the Douglas-Rachford method to distributedly find a solution. The proposed algorithm requires milder conditions compared to the existing methods. We prove the convergence of the proposed algorithm to an exact variational generalized Nash equilibrium under two different sets of assumptions. Our algorithm is validated numerically through the example of a Nash-Cournot production game.**

## I. Introduction

Recent years have witnessed increasing efforts in generalized Nash equilibrium problems (GNEP) [2], [3], motivated by numerous applications, e.g., communication networks [4], charge scheduling of electric vehicles [5], formation control [6], and demand management in the smart grid [7]. In many cases, multiple self-interested players/decision-makers aim to optimize their individual objectives under some global resource limits while unwilling to share their private information with the public. The goal is to design for each player a protocol that only requires information exchange within its trusted partners. By following it, these groups of players can eventually achieve a generalized Nash equilibrium (GNE).

Designing distributed solution for the GNEP has attracted lots of research interest. When players' local objectives only rely on the decisions of themselves and their trusted partners/neighbors, [8] proposes some elegant methods based on monotone operator theory. The work in [9], [10] extend these results by combining operator-splitting approaches with doubly-augmented information space. More precisely, these works build their approaches based on the "preconditioned" forward-backward method. To guarantee the convergence of the forward-backward algorithm, it is assumed that each player's objective is continuously differentiable in its decision given those of the others, and the associated pseudogradient/game Jacobian operator is strongly monotone and Lipschitz continuous. To relax some of these conditions, [11] proposes two algorithms based on alternating direction method of multipliers for a network with locally dependent objectives. The method in [12] solves locally linear coupled games using a proximal-point algorithm. [13] proposes a semi-decentralized method leveraging the Douglas-Rachford (DR) splitting method to handle aggregative games.

This paper focuses on the GNEP with generic objectives, which can be globally dependent and contain no special structure. The proposed algorithm only needs communication between each player and its trusted neighbors on a directed and weakly connected communication graph. Moreover, we establish the convergence of the proposed algorithm under two different assumptions on the pseudogradient and extended pseudogradient, allowing more flexibility compared with existing works. One assumption allows the local objective to be non-differentiable. Finally, we demonstrate the effectiveness of our proposed algorithm via a network Nash Cournot game.

*Basic Notations:* For a set of matrices $\{V_i\}_{i\in S}$, we let $\mathrm{blkd}(V_1,\dots,V_{|S|})$ denote the diagonal concatenation of these matrices, $[V_1,\dots,V_{|S|}]$ their horizontal stack, and $[V_1;\cdots;V_{|S|}]$ their vertical stack. For a set of vectors $\{v_i\}_{i\in S}$, $[v_i]_{i\in S}$ or $[v_1;\cdots;v_{|S|}]$ denotes their vertical stack. For a vector $v$ and a positive integer $i$, $[v]_i$ denotes the $i$th entry of $v$. Denote $\overline{\mathbb{R}} := \mathbb{R}\cup\{+\infty\}$, $\mathbb{R}_+ := [0,+\infty)$, and $\mathbb{R}_{++} := (0,+\infty)$. $\mathbb{S}^n_+$ (resp. $\mathbb{S}^n_{++}$) represents the set of all $n\times n$ symmetric positive semi-definite (resp. definite) matrices. $N_S(x)$ denotes the normal cone to the set $S\subseteq\mathbb{R}^n$ at the point $x$: if $x\in S$, then $N_S(x) := \{u\in\mathbb{R}^n \mid \sup_{z\in S}\langle u, z-x\rangle \le 0\}$; otherwise, $N_S(x) := \varnothing$. We use $\rightrightarrows$ to indicate a point-to-set map. For an operator $T:\mathbb{R}^n\rightrightarrows\mathbb{R}^n$, $\mathrm{Zer}(T) := \{x\in\mathbb{R}^n \mid Tx\ni\mathbf{0}\}$ and $\mathrm{Fix}(T) := \{x\in\mathbb{R}^n \mid Tx\ni x\}$ denote its zero set and fixed point set, respectively.

## II. Problem Formulation

### A. Game Formulation and GNE

We consider a set of players indexed by $\mathcal{N} := \{1,\dots,N\}$, each of which decides on its decision variables $x_i\in\Omega_i$ and optimizes its local objective function $J_i(x_i;x_{-i})$. Here, $\Omega_i\subseteq\mathbb{R}^{n_i}$ denotes player $i$'s local feasible set. Let the index set of players except player $i$ be defined by $\mathcal{N}_{-i}$. The vector $x_{-i}\in\mathbb{R}^{n_{-i}}$ represents the vertical stack of other players' decision variables, with $n_{-i} := \sum_{j\in\mathcal{N}_{-i}} n_j$. Formally, given the decision of other players, the player $i$ aims to solve a local optimization problem as follows:

$$\begin{cases}\underset{x_i\in\Omega_i}{\text{minimize}} & J_i(x_i;x_{-i})\\ \text{subject to} & A_ix_i \le b - \sum_{j\in\mathcal{N}_{-i}} A_jx_j\end{cases}, \tag{1}$$

This work was supported by the National Science Foundation under Grant No. 2014816. A conference version will be submitted to CDC 2021 [1].

The authors are with the School of Electrical and Computer Engineering, Purdue University, West Lafayette, IN, 47907, USA (e-mail: huan1282@purdue.edu; jianghai@purdue.edu).

where $A_i \in \mathbb{R}^{m \times n_i}$, $m$ is the number of the (global) affine coupling constraints, and $b \in \mathbb{R}^m$ is a constant vector. Denote the vertical stack of all decision variables by $x := [x_1; \cdots; x_N]$. The feasible set of the collective decision vector $x$, is given by:

$$\mathcal{X} := \Omega \cap \{x \in \mathbb{R}^n | Ax - b \leq \mathbf{0}\}, \tag{2}$$

where $\Omega := \prod_{i=1}^N \Omega_i$, $A := [A_1, A_2, \ldots, A_N]$, and $n = \sum_{i \in \mathcal{N}} n_i$. The feasible decision set of each player $i \in \mathcal{N}$ is characterized by the set-valued mapping $\mathcal{X}_i : \mathbb{R}^{n_{-i}} \rightrightarrows \mathbb{R}^{n_i}$, which is defined as:

$$\mathcal{X}_i(x_{-i}) := \{x_i \in \Omega_i | A_i x_i \leq b - \textstyle\sum_{j \in \mathcal{N}_{-i}} A_j x_j\}. \tag{3}$$

**Assumption 1.** *(Existence of Subgradient) For each $i \in \mathcal{N}$, $J_i(x_i; x_{-i})$ is a proper, lower semicontinuous, and convex function w.r.t. $x_i$, given any fixed $x_{-i}$.*

**Assumption 2.** *(Feasible Sets) Each local feasible set $\Omega_i$ is nonempty, closed, and convex. The global feasible set $\mathcal{X}$ is nonempty and satisfies Slater's constraint qualification [14, Ch 3.2].*

A generalized Nash equilibrium (GNE) for the game (1) is a joint decision vector $x$ such that, for each player $i$, $x_i \in \mathcal{X}_i(x_{-i})$ is a minimizer of the local optimization problem in (1). Equivalently, a GNE $x$ is a solution of the Karust-Kuhn-Tucker (KKT) conditions given as follows:

$$\begin{aligned} &\mathbf{0} \in \partial_{x_i} J_i(x_i; x_{-i}) + A_i^T \lambda_i + N_{\Omega_i}(x_i) \\ &\mathbf{0} \in -(Ax - b) + N_{\mathbb{R}_+^m}(\lambda_i), \end{aligned} \tag{4}$$

where $\lambda_i$ is the Lagrangian multiplier for the inequality constraints for the local problem (1) of player $i$. In this paper, we focus on the variational generalized Nash Equilibria (v-GNEs), which are a subset of GNEs [3]. The v-GNEs are solutions to the inclusions in (4) with all $\{\lambda_i\}_{i \in \mathcal{N}}$ on consensus ($\lambda_1 = \cdots = \lambda_N$). The KKT conditions for a v-GNE are:

$$\begin{aligned} &\mathbf{0} \in \partial_{x_i} J_i(x_i; x_{-i}) + A_i^T \lambda + N_{\Omega_i}(x_i) \\ &\mathbf{0} \in -(Ax - b) + N_{\mathbb{R}_+^m}(\lambda). \end{aligned} \tag{5}$$

for all $i \in \mathcal{N}$. Another way to characterize v-GNE is via variational inequalities (VI) [3]. Define the pseudogradient/game Jacobian $F : \mathbb{R}^n \rightrightarrows \mathbb{R}^n$ as:

$$F : x \mapsto [\partial_{x_i} J_i(x_i; x_{-i})]_{i \in \mathcal{N}}. \tag{6}$$

An associated variational inequality problem VI$(\mathcal{X}, F)$ can accordingly be defined as: finding $x^* \in \mathcal{X}$ such that $\langle F(x^*), x - x^* \rangle \geq 0, \forall x \in \mathcal{X}$ (see [3]). Solutions to VI$(\mathcal{X}, F)$ are exactly the v-GNE of problem (1), whose existence and uniqueness have been well studied [14].

**Remark 1.** *We focus on v-GNE since we can leverage the rich body of theory and tools developed for solving VIs [14] [15, Ch 12], and by keeping all $\{\lambda_i\}_{i \in \mathcal{N}}$ on consensus, v-GNE possesses desirable properties such as economic fairness and better social stability/sensitivity [16].*

### *B. Networked Game Formulation*

To enable the distributed computation of v-GNE, we consider an underlying communication graph $\mathcal{G} = (\mathcal{N}_g, \mathcal{E}_g)$, where players can communicate with their neighbors through arbitrators on the directed edges. The node set $\mathcal{N}_g$ represents the set of all the players, and $\mathcal{E}_g \subseteq \mathcal{N}_g \times \mathcal{N}_g$ is the set of directed edges. The cardinalities $|\mathcal{N}_g|$ and $|\mathcal{E}_g|$ are denoted by $N_g$ and $E_g$. In this case, $\mathcal{N}_g = \mathcal{N}$ and $N_g = N$. We use $(i, j)$ to denote a directed edge having node/player $i$ as its tail and node/player $j$ as its head.

**Assumption 3.** *(Communicability) The underlying communication graph $\mathcal{G} = (\mathcal{N}_g, \mathcal{E}_g)$ is directed and weakly connected. Besides, it has no self-loops.*

Let $\mathcal{N}_i$ denote the set of immediate neighbors of player $i$ who can directly communicate with it, $\mathcal{N}_i^+ := \{j \in \mathcal{N} \mid (j, i) \in \mathcal{E}_g\}$ the set of in-neighbors of player $i$, and $\mathcal{N}_i^- := \{j \in \mathcal{N} \mid (i, j) \in \mathcal{E}_g\}$ the set of out-neighbors of player $i$. Although the graph $\mathcal{G}$ is directed, we assume that each node can send messages to both its in- and out-neighbors.

Our goal is to convert the centralized GNEP in (1) into the zero-finding problem of a certain operator that can be carried out distributedly over the graph $\mathcal{G}$. To construct this operator, in the rest of this section, we will first derive the network Lagrangian. Let each player extend its decision space to maintain the local estimates of other players' decision vectors. Each player $i$ now has augmented decision variable $y_i \in \mathbb{R}^n$, which consists of its local decision $y_i^i \in \mathbb{R}^{n_i}$ and the local estimate $y_i^j \in \mathbb{R}^{n_j}$ of player $j$'s decision for each $j \in \mathcal{N}_{-i}$. Let $y_i^{-i} \in \mathbb{R}^{n_{-i}}$ denote the vertical stack of $y_i^j$ for all $j \in \mathcal{N}_{-i}$ in a prespecified order. Denote $n_{<i} = \sum_{j \in \mathcal{N}, j < i} n_j$ and $n_{>i} = \sum_{j \in \mathcal{N}, j > i} n_j$. The extended feasible set $\tilde{\Omega}$ is defined as $\tilde{\Omega} := \tilde{\Omega}_1 \times \tilde{\Omega}_2 \times \cdots \times \tilde{\Omega}_N$, where the feasible set of each $y_i$ is defined as $\tilde{\Omega}_i := \mathbb{R}^{n_{<i}} \times \Omega_i \times \mathbb{R}^{n_{>i}}$. Notice that $\mathbb{R}^{n_{<i}} \times N_{\Omega_i}(y_i^i) \times \mathbb{R}^{n_{>i}} = N_{\tilde{\Omega}_i}(y_i)$.

For the purpose of deducing the network Lagrangian only (but not for eventual distributed algorithm development), we select an arbitrary player $\nu$ to be the conceptual leader to take care of the global resource constraints, whose local problem is given by:

$$\begin{cases} \underset{y_\nu \in \tilde{\Omega}_\nu}{\text{minimize}} & J_\nu(y_\nu^\nu; u_\nu)_{|u_\nu = y_\nu^{-\nu}} \\ \text{subject to} & A_\nu y_\nu^\nu \leq b - \sum_{j \in \mathcal{N}_{-\nu}} A_j y_\nu^j \\ & y_\nu = y_j, \forall j \in \mathcal{N}_\nu^+ \end{cases}, \tag{7}$$

**Remark 2.** *Notice that in the original game formulation (1), $x_{-i}$ is treated as a parametric input in player $i$'s local optimization problem. In (7), other players' decision vectors $x_{-\nu}$ are replaced with $y_\nu^{-\nu}$; hence $y_\nu^{-\nu}$ is also treated as a parameter in the objective function $J_\nu$.*

With player $\nu$ taking care of the global resource constraints, the remaining players $i \in \mathcal{N}_{-\nu}$ only need to cope with optimization problems with consensus constraints:

$$\begin{cases} \underset{y_i \in \tilde{\Omega}_i}{\text{minimize}} & J_i(y_i^i; u_i)_{|u_i = y_i^{-i}} \\ \text{subject to} & y_i = y_j, \forall j \in \mathcal{N}_i^+ \end{cases}. \tag{8}$$

The Lagrangian for the interdependent constrained optimization problem of player $\nu$ is given by:

$$\mathcal{L}_\nu(y_\nu, \{y_j\}_{j\in\mathcal{N}_{-\nu}}, \lambda, \{\mu_{j\nu}\}_{j\in\mathcal{N}_{-\nu}}) := J_\nu(y_\nu^\nu; u_\nu)_{|u_\nu=y_\nu^{-\nu}} + N_{\Omega_\nu}(y_\nu^\nu) + \lambda^T(Ay_\nu - b) + \textstyle\sum_{j\in\mathcal{N}_\nu^+}\mu_{j\nu}^T(y_\nu - y_j). \tag{9}$$

The Lagrangian of player $i$, where $i \in \mathcal{N}_{-\nu}$, is

$$\mathcal{L}_i(y_i, \{y_j\}_{j\in\mathcal{N}_{-i}}, \{\mu_{ji}\}_{j\in\mathcal{N}_{-i}}) := J_i(y_i^i; u_i)_{|u_i=y_i^{-i}} + N_{\Omega_i}(y_i^i) + \textstyle\sum_{j\in\mathcal{N}_i^+}\mu_{ji}^T(y_i - y_j). \tag{10}$$

By summing these individual Lagrangians, we obtain the network Lagrangian as follows:

$$\mathcal{L}_{\text{net}}(\{y_i\}_{i\in\mathcal{N}}, \{\mu_{ji}\}_{(j,i)\in\mathcal{E}_g}, \lambda) := \lambda^T(Ay_\nu - b) + \sum_{i\in\mathcal{N}}\left(J_i(y_i^i; u_i)_{u_i=y_i^{-i}} + N_{\Omega_i}(y_i^i)\right) + \sum_{(j,i)\in\mathcal{E}_g}\mu_{ji}^T(y_i - y_j). \tag{11}$$

The above proposed scheme has the drawbacks that the leader $\nu$ is assumed to know $A_i$ from every other player $i$, and is the only one in charge of the global resource constraint (hence prune to single-point failures). To avoid these issues, we make the following modifications to $\mathcal{L}_{\text{net}}$. First, let each player $i$ keep a local estimate $\lambda_i \in \mathbb{R}_+^m$ of $\lambda$, along with additional constraints to guarantee $\{\lambda_i\}_{i\in\mathcal{N}}$ are on consensus. Second, decompose the global constraint $Ay_\nu - b \le \mathbf{0}$ into $\sum_i(A_iy_i^i - b_i) \le \mathbf{0}$, where each $y_\nu^i$ is replaced with $y_i^i$, and $\{b_i\}_{i\in\mathcal{N}} \in \mathbb{R}^m$ is a set of arbitrary vectors satisfying $\sum_{i\in\mathcal{N}} b_i = b$. Also, to facilitate the convergence, we introduce some second order penalty terms for consensus constraints. The modified augmented Lagrangian is given as follows:

$$\begin{aligned}\mathcal{L}_{\text{net}}^* := &\sum_{i\in\mathcal{N}}\left(J_i(y_i^i; u_i)_{|u_i=y_i^{-i}} + N_{\Omega_i}(y_i^i) + \lambda_i^T(A_iy_i^i - b_i)\right)\\ &+ \sum_{(j,i)\in\mathcal{E}_g}\left(\mu_{ji}^T(y_i - y_j) - z_{ji}^T(\lambda_i - \lambda_j)\right)\\ &+ \sum_{(j,i)\in\mathcal{E}_g}\left(\frac{\rho_\mu}{2}\|y_i - y_j\|^2 - \frac{\rho_z}{2}\|\lambda_i - \lambda_j\|^2\right),\end{aligned} \tag{12}$$

where $\{z_{ji}\}_{(j,i)\in\mathcal{E}_g}$ is the set of Lagrange multipliers ensuring the consensus of $\{\lambda_i\}_{i\in\mathcal{N}}$. We let two constant parameters $\rho_\mu$ and $\rho_z$ control the weights of the second order penalty terms, with lower bounds to be determined later on.

In the augmented network Lagrangian, each player $i$'s objective function $J_i(y_i^i; y_i^{-i})$ is only optimized over $y_i^i$, while its local estimates $y_i^{-i}$ of other players are updated to satisfy the global and consensus constraints. Accordingly, the extended pseudogradient $\mathcal{F} : \mathbb{R}^{nN} \rightrightarrows \mathbb{R}^n$ is the set-valued operator defined as

$$\mathcal{F} : \boldsymbol{y} \mapsto [\partial_{y_i^i} J_i(y_i^i; y_i^{-i})]_{i\in\mathcal{N}}, \tag{13}$$

where $\boldsymbol{y} := [y_1; \cdots; y_N] \in \mathbb{R}^{nN}$. For $\boldsymbol{y}$ having $\{y_i\}_{i\in\mathcal{N}}$ on consensus, $F(y_i) = \mathcal{F}(\boldsymbol{y})$. To incorporate the extended pseudogradient $\mathcal{F}$ into a fixed point iteration, we introduce the individual selection matrices $\{\mathcal{R}_i\}_{i\in\mathcal{N}}$ and their diagonal concatenation $\mathcal{R} \in \mathbb{R}^{n\times nN}$:

$$\mathcal{R}_i = [\mathbf{0}_{n_i\times n_{<i}}, \mathbf{I}_{n_i}, \mathbf{0}_{n_i\times n_{>i}}], \mathcal{R} = \text{blkd}(\mathcal{R}_1, \dots, \mathcal{R}_N). \tag{14}$$

Notice that $y_i^i = \mathcal{R}_iy_i$ and $\mathcal{R}_i\mathcal{R}_i^T = I_{n_i}$.

## III. Distributed Algorithm with the DR Splitting

In this section, we recast the GNEP to a zero-finding problem for a properly defined operator. We derive the analytical updating steps to solve for the zeros by leveraging the Douglas-Rachford operator splitting method and constructing a design matrix. The analysis of monotonicity and convergence will be discussed in Section IV.

### A. Zero-finding Problem

To solve the distributed GNE seeking problem, we need to find the stationary points of the augmented network Lagrangian (12). By taking the partial derivative of $\mathcal{L}_{\text{net}}^*$ w.r.t. each variable and reversing the sign of the rows corresponding to $\lambda_i$ and $\mu_i$, we can obtain the following set-valued operator $\mathcal{T}$:

$$\mathcal{T} : \begin{bmatrix}\boldsymbol{y}\\ \boldsymbol{\lambda}\\ \boldsymbol{\mu}\\ \boldsymbol{z}\end{bmatrix} \mapsto \begin{bmatrix}\mathcal{R}^T(\mathcal{F}(\boldsymbol{y}) + \Lambda^T\boldsymbol{\lambda}) + B_n\boldsymbol{\mu} + \rho_\mu L_n\boldsymbol{y} + N_{\tilde{\Omega}}(\boldsymbol{y})\\ N_{\mathbb{R}_+^{mN}}(\boldsymbol{\lambda}) - \Lambda\mathcal{R}\boldsymbol{y} + \boldsymbol{b} + B_m\boldsymbol{z} + \rho_z L_m\boldsymbol{\lambda}\\ -B_n^T\cdot\boldsymbol{y}\\ -B_m^T\cdot\boldsymbol{\lambda}\end{bmatrix}, \tag{15}$$

where $\Lambda$ is the diagonal concatenation of $\{A_i\}_{i\in\mathcal{N}}$, i.e., $\Lambda := \text{blkd}(A_1, \dots, A_N)$; $\boldsymbol{b}$ is the vertical stack of $\{b_i\}_{i\in\mathcal{N}}$; $B_n := (B\otimes I_n)$, $L_n := (L\otimes I_n)$, $B_m := (B\otimes I_m)$, $L_m := (L\otimes I_m)$, $B$ and $L$ are the incidence matrix and Laplacian matrix of the underlying communication graph, respectively, with $L = B\cdot B^T$; and $\boldsymbol{\lambda}$, $\boldsymbol{\mu}$, and $\boldsymbol{z}$ are the stack vectors of $\{\lambda_i\}_{i\in\mathcal{N}}$, $\{\mu_{ji}\}_{(j,i)\in\mathcal{E}_g}$, and $\{z_{ji}\}_{(j,i)\in\mathcal{E}_g}$, respectively.

**Theorem 1.** *Suppose Assumptions 1-3 hold, and there exists $\omega^* = [\boldsymbol{y}^*; \boldsymbol{\lambda}^*; \boldsymbol{\mu}^*; \boldsymbol{z}^*] \in Zer(\mathcal{T})$. Then $\boldsymbol{y}^* = \mathbf{1}_N \otimes y^*$, $\boldsymbol{\lambda}^* = \mathbf{1}_N \otimes \lambda^*$, and $(y^*, \lambda^*)$ satisfies the KKT conditions* (5) *for v-GNE with x replaced with $y^*$. Furthermore, for a v-GNE $(\boldsymbol{y}^\dagger, \boldsymbol{\lambda}^\dagger)$ of the game* (1)*, there exist $\boldsymbol{\mu}^\dagger$ and $\boldsymbol{z}^\dagger$ such that $[\boldsymbol{y}^\dagger; \boldsymbol{\lambda}^\dagger; \boldsymbol{\mu}^\dagger; \boldsymbol{z}^\dagger] \in Zer(\mathcal{T})$.*

*Proof.* See Appendix A. □

Theorem 1 implies that we can convert the solution of v-GNE of the original GNEP into the zero-finding problem of the operator $\mathcal{T}$ defined in (15). In the next subsection, we will propose a candidate algorithm to solve the latter problem.

### B. Operator Splitting

Given a generic set-valued operator $T : \mathbb{R}^\ell \rightrightarrows \mathbb{R}^\ell$, $T$ is monotone if, for any $(x, u) \in \text{gra}(T)$ and $(x', u') \in \text{gra}(T)$, $\langle x - x', u - u'\rangle \ge 0$, where $\text{gra}(T)$ is the graph of the operator $T$. If $T$ is maximally monotone, a point in $\text{Zer}(T)$ could in principle be determined through the proximal-point algorithm, namely, the fixed point iteration using its resolvent operator $J_T : \mathbb{R}^\ell \to \mathbb{R}^\ell$ defined as $J_T := (I + T)^{-1}$ [17, Thm. 23.41]. However, it is often infeasible to evaluate the resolvent in a distributed manner for operators arising in network problems such as the one in (15). For this purpose, we consider the Douglas-Rachford (DR) splitting technique, where $T = A + B$ is split into two maximally monotone operators $A$ and $B$ whose resolvents can be more conveniently computed. The resolvent and reflected resolvent of $A$ are defined as $J_A := (I + A)^{-1}$ and $R_A := 2J_A - I$, respectively. As $A$ is maximally monotone, $J_A$ is firmly nonexpansive and $R_A$ is non-expansive, i.e., for any $x$ and $x'$, $\|J_Ax-J_Ax'\|^2+\|(I-J_A)x-(I-J_A)x'\|^2 \le \|x-x'\|^2$ and $\|R_Ax-R_Ax'\| \le \|x-x'\|^2$ (see [17, Cor. 23.11]). Similar properties hold for $B$.

**Algorithm 1:** Douglas-Rachford Splitting

**Initialize:** $\tilde{\omega}^{(0)}$;
**Iterate until convergence:**
$\omega^{(k+1)} = J_{\tilde{\mathcal{A}}}(\tilde{\omega}^{(k)})$ ;
$\tilde{\omega}^{(k+1)} = \tilde{\omega}^{(k)} + \Big(J_{\tilde{\mathcal{B}}}(2\omega^{(k+1)} - \tilde{\omega}^{(k)}) - \omega^{(k+1)}\Big)$;
**Return:** $\omega^{(k)}$.

The DR splitting algorithm can be viewed as a special case of the Krasnosel'skii-Mann algorithm. Given a nonexpansive operator $Q$ which has a nonempty fixed point set $\mathrm{Fix}(Q)$, the Krasnosel'skii-Mann algorithm suggests the iteration:

$$x^{(k+1)} = x^{(k)} + \gamma^{(k)}(Qx^{(k)} - x^{(k)}). \tag{16}$$

If $\{\gamma^{(k)}\}_{k\in\mathbb{N}}$ is a sequence in $[0,1]$ satisfying $\sum_{k\in\mathbb{N}}\gamma^{(k)}(1-\gamma^{(k)}) = \infty$, then $\{x^{(k)}\}_{k\in\mathbb{N}}$ converges to a point in $\mathrm{Fix}(Q)$ [17, Ch. 5.2]. In the DR algorithm, $Q$ in (16) is set to be $R_BR_A$. Then, the iteration in (16) converges to some $x^* \in \mathrm{Fix}(R_BR_A)$. It follows from [17, Ch. 26.3] that $J_A(x^*) \in \mathrm{Zer}(A+B) = \mathrm{Zer}(T)$.

Now we focus on the operator $\mathcal{T}$ in (15) and split it into the operators $\mathcal{A}$ and $\mathcal{B}$ defined as follows:

$$\mathcal{A}: \begin{bmatrix} \boldsymbol{y} \\ \boldsymbol{\lambda} \\ \boldsymbol{\mu} \\ \boldsymbol{z} \end{bmatrix} \mapsto \mathcal{A}_y \begin{bmatrix} \boldsymbol{y} \\ \boldsymbol{\lambda} \\ \boldsymbol{\mu} \\ \boldsymbol{z} \end{bmatrix} + D \begin{bmatrix} \boldsymbol{y} \\ \boldsymbol{\lambda} \\ \boldsymbol{\mu} \\ \boldsymbol{z} \end{bmatrix} + \begin{bmatrix} N_{\tilde{\Omega}}(\boldsymbol{y}) \\ N_{\mathbb{R}_+^{mN}}(\boldsymbol{\lambda}) \\ 0 \\ 0 \end{bmatrix} \tag{17}$$

and

$$\mathcal{B}: \begin{bmatrix} \boldsymbol{y} \\ \boldsymbol{\lambda} \\ \boldsymbol{\mu} \\ \boldsymbol{z} \end{bmatrix} \mapsto D \begin{bmatrix} \boldsymbol{y} \\ \boldsymbol{\lambda} \\ \boldsymbol{\mu} \\ \boldsymbol{z} \end{bmatrix} + \begin{bmatrix} \frac{\rho_\mu}{2}L_n & 0 & 0 & 0 \\ 0 & \frac{\rho_z}{2}L_m & 0 & 0 \\ 0 & 0 & 0 & 0 \\ 0 & 0 & 0 & 0 \end{bmatrix} \begin{bmatrix} \boldsymbol{y} \\ \boldsymbol{\lambda} \\ \boldsymbol{\mu} \\ \boldsymbol{z} \end{bmatrix}, \tag{18}$$

with

$$D = \begin{bmatrix} 0 & \frac{1}{2}(\Lambda\mathcal{R})^T & \frac{1}{2}B_n & 0 \\ -\frac{1}{2}\Lambda\mathcal{R} & 0 & 0 & \frac{1}{2}B_m \\ -\frac{1}{2}B_n^T & 0 & 0 & 0 \\ 0 & -\frac{1}{2}B_m^T & 0 & 0 \end{bmatrix}, \tag{19}$$

and

$$\mathcal{A}_y: \begin{bmatrix} \boldsymbol{y} \\ \boldsymbol{\lambda} \\ \boldsymbol{\mu} \\ \boldsymbol{z} \end{bmatrix} \mapsto \begin{bmatrix} \mathcal{R}^T\mathcal{F}(\boldsymbol{y}) + \frac{\rho_\mu}{2}L_n\boldsymbol{y} \\ \boldsymbol{b} + \frac{\rho_z}{2}L_m\boldsymbol{\lambda} \\ 0 \\ 0 \end{bmatrix}. \tag{20}$$

Given the splitted operators $\mathcal{A}$ and $\mathcal{B}$, to allow for the distributed evaluation of the resolvents $J_{\mathcal{A}}$ and $J_{\mathcal{B}}$, a positive definite matrix $\Phi$ is designed. It is constructed in a way such that $D+\Phi$ is lower or upper triangular. One choice is

$$\Phi = \begin{bmatrix} \boldsymbol{\tau}_1^{-1} - \frac{\rho_\mu}{2}L_n & -\frac{1}{2}(\Lambda\mathcal{R})^T & -\frac{1}{2}B_n & 0 \\ -\frac{1}{2}\Lambda\mathcal{R} & \boldsymbol{\tau}_2^{-1} - \frac{\rho_z}{2}L_m & 0 & -\frac{1}{2}B_m \\ -\frac{1}{2}B_n^T & 0 & \boldsymbol{\tau}_3^{-1} & 0 \\ 0 & -\frac{1}{2}B_m^T & 0 & \boldsymbol{\tau}_4^{-1} \end{bmatrix}, \tag{21}$$

where $\boldsymbol{\tau}_1 := \mathrm{blkd}(\tau_{11}\otimes I_n, \ldots, \tau_{1N}\otimes I_n)$; similarly for $\boldsymbol{\tau}_2$, $\boldsymbol{\tau}_3$ and $\boldsymbol{\tau}_4$.

Let $\tilde{\mathcal{A}}$ denote $\Phi^{-1}\mathcal{A}$ and $\tilde{\mathcal{B}}$ denote $\Phi^{-1}\mathcal{B}$. Note that $\mathbf{0} \in (\mathcal{A}+\mathcal{B})\omega$ if and only if $\mathbf{0} \in (\tilde{\mathcal{A}}+\tilde{\mathcal{B}})\omega$. By applying the iteration in (16)) to $\tilde{\mathcal{A}}$ and $\tilde{\mathcal{B}}$ (let $Q = R_{\tilde{\mathcal{B}}}R_{\tilde{\mathcal{A}}}$) and setting $\gamma_k$ always equal to $\frac{1}{2}$, we can obtain Algorithm 1. A detailed version with explicit operations at each node and edge is given in Algorithm 2. The result $\{y_i^{(k)}\}_{i\in\mathcal{N}}$ returned by Algorithm 2 will be proved in the next section to converge to a v-GNE of the problem (1).

**Algorithm 2:** Distributed GNE Seeking

**Initialize:** $\{\tilde{y}_i^{(0)}\}, \{\tilde{\lambda}_i^{(0)}\}, \{\tilde{\mu}_{ji}^{(0)}\}, \{\tilde{z}_{ji}^{(0)}\}$;
**Iterate until convergence:**
**for** *player* $i \in \mathcal{N}$ **do**

$$y_i^{-i(k+1)} = \tilde{y}_i^{-i(k)} - \tau_{1i}\big(\tfrac{\rho_\mu}{2}\textstyle\sum_{j\in\mathcal{N}_i}\big(\tilde{y}_i^{-i(k)} - \tilde{y}_j^{-i(k)}\big) + \tfrac{1}{2}\big(\sum_{j\in\mathcal{N}_i^+}\tilde{\mu}_{ji}^{-i(k)} - \sum_{j\in\mathcal{N}_i^-}\tilde{\mu}_{ij}^{-i(k)}\big)\big) ;$$

$$y_i^{i(k+1)} = \underset{v_i\in\Omega_i}{\mathrm{argmin}}\Big[J_i\big(v_i; y_i^{-i(k+1)}\big) + \tfrac{1}{2}(\tilde{\lambda}_i^{(k)})^T A_i v_i + \tfrac{\rho_\mu}{2}\sum_{j\in\mathcal{N}_i}\big(\tilde{y}_i^{i(k)} - \tilde{y}_j^{i(k)}\big)^T v_i + \tfrac{1}{2\tau_{1i}}\big\|v_i - \tilde{y}_i^{i(k)}\big\|^2 + \tfrac{1}{2}\big(\textstyle\sum_{j\in\mathcal{N}_i^+}\tilde{\mu}_{ji}^{i(k)} - \sum_{j\in\mathcal{N}_i^-}\tilde{\mu}_{ij}^{i(k)}\big)^T v_i\Big];$$

$$\lambda_i^{(k+1)} = \mathrm{Proj}_{\mathbb{R}_+^m}\Big[\tilde{\lambda}_i^{(k)} + \tau_{2i}\big(A_i\mathcal{R}_i\cdot y_i^{(k+1)} - \tfrac{1}{2}A_i\mathcal{R}_i\tilde{y}_i^{(k)} - \tfrac{\rho_z}{2}\textstyle\sum_{j\in\mathcal{N}_i}(\tilde{\lambda}_i^{(k)} - \tilde{\lambda}_j^{(k)}) - \tfrac{1}{2}(\sum_{j\in\mathcal{N}_i^+}\tilde{z}_{ji}^{(k)} - \sum_{j\in\mathcal{N}_i^-}\tilde{z}_{ij}^{(k)}) - b_i\big)\Big];$$

**end**
**for** *edge* $(j,i) \in \mathcal{E}_g$ **do**

$$\mu_{ji}^{(k+1)} = \tilde{\mu}_{ji}^{(k)} + \tau_{3i}\big(y_i^{(k+1)} - y_j^{(k+1)} - \tfrac{1}{2}(\tilde{y}_i^{(k)} - \tilde{y}_j^{(k)})\big);$$

$$z_{ji}^{(k+1)} = \tilde{z}_{ji}^{(k)} + \tau_{4i}\big(\lambda_i^{(k+1)} - \lambda_j^{(k+1)} - \tfrac{1}{2}(\tilde{\lambda}_i^{(k)} - \tilde{\lambda}_j^{(k)})\big);$$

**end**
**for** *player* $i \in \mathcal{N}$ **do**

$$\tilde{y}_i^{(k+1)} = y_i^{(k+1)} - \tau_{1i}\big(\mathcal{R}_i^T A_i^T\big(\lambda_i^{(k+1)} - \tfrac{1}{2}\tilde{\lambda}_i^{(k)}\big) + \rho_\mu\textstyle\sum_{j\in\mathcal{N}_i}\big((y_i^{(k+1)} - y_j^{(k+1)}) - \tfrac{1}{2}(\tilde{y}_i^{(k)} - \tilde{y}_j^{(k)})\big) + \sum_{j\in\mathcal{N}_i^+}(\mu_{ji}^{(k+1)} - \tfrac{1}{2}\tilde{\mu}_{ji}^{(k)}) - \sum_{j\in\mathcal{N}_i^-}(\mu_{ij}^{(k+1)} - \tfrac{1}{2}\tilde{\mu}_{ij}^{(k)})\big);$$

$$\tilde{\lambda}_i^{(k+1)} = \lambda_i^{(k+1)} + \tau_{2i}\big(A_i\mathcal{R}_i\big(\tilde{y}_i^{(k+1)} - \tfrac{1}{2}\tilde{y}_i^{(k)}\big) - \rho_z\textstyle\sum_{j\in\mathcal{N}_i}\big((\lambda_i^{(k+1)} - \lambda_j^{(k+1)}) - \tfrac{1}{2}(\tilde{\lambda}_i^{(k)} - \tilde{\lambda}_j^{(k)})\big) - \sum_{j\in\mathcal{N}_i^+}(z_{ji}^{(k+1)} - \tfrac{1}{2}\tilde{z}_{ji}^{(k)}) + \sum_{j\in\mathcal{N}_i^-}(z_{ij}^{(k+1)} - \tfrac{1}{2}\tilde{z}_{ij}^{(k)})\big);$$

**end**
**for** *edge* $(j,i) \in \mathcal{E}_g$ **do**

$$\tilde{\mu}_i^{(k+1)} = \mu_i^{(k+1)} + \tau_{3i}\big(\tilde{y}_i^{(k+1)} - \tilde{y}_j^{(k+1)} - \tfrac{1}{2}(\tilde{y}_i^{(k)} - \tilde{y}_j^{(k)})\big);$$

$$\tilde{z}_i^{(k+1)} = z_i^{(k+1)} + \tau_{4i}\big(\tilde{\lambda}_i^{(k+1)} - \tilde{\lambda}_j^{(k+1)} - \tfrac{1}{2}(\tilde{\lambda}_i^{(k)} - \tilde{\lambda}_j^{(k)})\big);$$

**end**
**Return:** $\{y_i^{(k)}\}_{i\in\mathcal{N}}$.

In Algorithm 2, the computational workload is concentrated on the update of local decision $y_i^i$ inside the first for-loop, while the other updates are computationally trivial. At each iteration, player $i$ will receive $(2n+2m)$ variables from each neighbor and send the same amount back to it. At the same time, player $i$ will also receive and send $(2n+2m)$ variables from and to the arbitrator on each incident edge.

## IV. Convergence Results

To analyze the convergence properties of the Algorithm 2, we make two parallel assumptions, either of which can guarantee the convergence to a v-GNE. Moreover, depending on the specific game scenario, either one could be less restrictive than the other one.

**Assumption 4.** *(Convergence Condition A) The GNEP admits at least one v-GNE, and the operator $\mathcal{R}^T\mathcal{F} + \frac{\rho_\mu}{2}L_n$ is maximally monotone.*

**Assumption 5.** *(Convergence Condition B) The pseudogradient operator $F$ is strongly monotone and Lipschitz continuous, i.e., there exist $\eta > 0$ and $\theta_1 > 0$, such that $\forall x, x' \in \mathbb{R}^n$, $\langle x - x', F(x) - F(x')\rangle \geq \eta\|x - x'\|^2$ and $\|F(x) - F(x')\| \leq \theta_1\|x - x'\|$. The operator $\mathcal{R}^T\mathcal{F}$ is Lipschitz continuous, i.e., there exists $\theta_2 > 0$, such that $\forall \boldsymbol{y}, \boldsymbol{y}' \in \mathbb{R}^{nN}$, $\|\mathcal{F}(\boldsymbol{y}) - \mathcal{F}(\boldsymbol{y}')\| \leq \theta_2\|\boldsymbol{y} - \boldsymbol{y}'\|$.*

**Remark 3.** *The Lipschitz continuity of $\mathcal{F}$ implies that of $F$. To see this, $\forall \boldsymbol{y}, \boldsymbol{y}' \in \mathbb{R}^{nN}$ and assuming their local estimates are on consensus ($y_i = y_j$ and $y_i' = y_j'$ for all $i, j$), we have $\|F(y_i) - F(y_i')\| = \|\mathcal{F}(\boldsymbol{y}) - \mathcal{F}(\boldsymbol{y}')\| \leq \theta_2\|\boldsymbol{y} - \boldsymbol{y}'\| = \sqrt{N}\cdot\theta_2\|y_i - y_i'\|$. However, in some cases, $\theta_1$ could provide a tighter Lipschitz constant than $\sqrt{N}\cdot\theta_2$, i.e., $\|F(y_i) - F(y_i')\| \leq \theta_1\|y_i - y_i'\| \leq \sqrt{N}\cdot\theta_2\|y_i - y_i'\|$.*

**Lemma 1.** *Suppose $\{\boldsymbol{\tau}_1, \boldsymbol{\tau}_2, \boldsymbol{\tau}_3, \boldsymbol{\tau}_4\}$ in the design matrix $\Phi$ satisfy the following inequalities: $\tau_{1i}^{-1} > \frac{1}{2}\|A_i\|_1 + (\frac{1}{2} + \rho_\mu)d_i$, $\tau_{2i}^{-1} > \frac{1}{2}\|A_i\|_\infty + (\frac{1}{2} + \rho_z)d_i, \forall i \in \mathcal{N}$, and $\tau_{3j}^{-1} > 1$, $\tau_{4j}^{-1} > 1$, $\forall j \in \mathcal{E}_g$. Then $\Phi$ is positive definite.*

Here, $d_i$ denotes the degree of node/player $i$. Lemma 1 is the direct result of the Gershgorin circle theorem: a sufficient condition for the design matrix to be positive definite is that all of its Gershgorin discs lie on the positive orthant [18]. The proof is trivial, and therefore omitted here.

Let $\mathcal{K}$ be the inner product space obtained by endowing the vector space $\mathbb{R}^{n(N+E_g)+m(N+E_g)}$ with the inner product $\langle\omega, \omega'\rangle_{\mathcal{K}} = \langle\Phi\omega, \omega'\rangle$.

**Lemma 2.** *Suppose the design matrix $\Phi$ is positive definite. Then the operator $\tilde{\mathcal{B}}$ is maximally monotone on $\mathcal{K}$, and its reflected resolvent $R_{\tilde{\mathcal{B}}}$ is nonexpansive on $\mathcal{K}$.*

*Proof.* $\mathcal{B}$ is a linear map whose matrix is the sum of the skew-symmetric matrix $D$ and another positive semi-definite matrix, and hence is maximally monotone [17, Ep. 20.35]. Since $\Phi$ is positive definite, it follows from [17, Prop. 20.24] that $\tilde{\mathcal{B}}$ is maximally monotone on $\mathcal{K}$. By [17, Cor. 23.11], this implies that the reflected resolvent $R_{\tilde{\mathcal{B}}} := 2J_{\tilde{\mathcal{B}}} - I$ is nonexpansive on $\mathcal{K}$. □

**Lemma 3.** *Suppose the design matrix $\Phi$ is positive definite and Assumptions 2 and 4 hold. Then the operator $\tilde{\mathcal{A}}$ is maximally monotone on $\mathcal{K}$, and its reflected resolvent $R_{\tilde{\mathcal{A}}}$ is nonexpansive on $\mathcal{K}$.*

*Proof.* By Assumption 2 and 4, $\mathcal{R}^T\mathcal{F} + \frac{\rho_\mu}{2}L_n$ and $N_{\tilde{\Omega}}$ are maximally monotone operators. The same arguments in the proof of Lemma 2 apply here. □

**Theorem 2.** *Suppose that Assumptions 1-4 hold and the design matrix $\Phi$ satisfies the inequalities in Lemma 1. Then the sequence $\{\boldsymbol{y}_{(k)}\}_{k\in\mathbb{N}}$ and $\{\boldsymbol{\lambda}_{(k)}\}_{k\in\mathbb{N}}$ generated by Algorithm 2 satisfy $\lim_{k\to\infty}\boldsymbol{y}_k = (\mathbf{1}_N \otimes y^*)$ and $\lim_{k\to\infty}\boldsymbol{\lambda}_k = (\mathbf{1}_N \otimes \lambda^*)$, where $y^*$ is a v-GNE to problem (1) and $(y^*, \lambda^*)$ together is a solution to the KKT conditions (5).*

Theorem 2 directly follows from the fact that $\tilde{\mathcal{A}}$ and $\tilde{\mathcal{B}}$ are two maximally monotone operators on $\mathcal{K}$ and the standard D-R results [17, Thm. 26.11]. Next, we are going to analyze the convergence properties under the Assumption 5 instead using the key notion of restricted monotonicity [10]. Given a generic operator $T$, we say that $T$ is restricted monotone w.r.t. a set $S$ if, for all $x$ and all $x^* \in S$, $\langle x - x^*, T(x) - T(x^*)\rangle \geq 0$. Restricted nonexpansiveness and restricted firm nonexpansiveness are defined similarly. The main results will be given in Theorem 3 below. First we introduce some preliminary results.

**Lemma 4.** *Suppose Assumptions 1-3 and 5 hold and $\rho_\mu \geq \frac{2}{\sigma_1}(\frac{(\theta_1+\theta_2)^2}{4\eta} + \theta_2)$. Then the operator $\mathcal{A}_y$ is restricted monotone w.r.t. Zer($\mathcal{T}$), i.e., for any $\omega$ and $\omega^*$ with $\omega^* \in$ Zer($\mathcal{T}$), $\langle\omega - \omega^*, \mathcal{A}_y(\omega) - \mathcal{A}_y(\omega^*)\rangle \geq 0$.*

*Proof.* See Appendix B. □

**Lemma 5.** *Suppose Assumptions 1-3 and 5 hold and $\rho_\mu \geq \frac{2}{\sigma_1}(\frac{(\theta_1+\theta_2)^2}{4\eta} + \theta_2)$. Then the operator $\tilde{\mathcal{A}}$ is restricted monotone w.r.t. Zer($\mathcal{T}$) on $\mathcal{K}$, and its resolvent $J_{\tilde{\mathcal{A}}}$ is restricted firmly nonexpansive w.r.t. Fix($R_{\tilde{\mathcal{B}}}R_{\tilde{\mathcal{A}}}$) on $\mathcal{K}$. Furthermore, the reflected resolvent $R_{\tilde{\mathcal{A}}}$ is restricted nonexpansive w.r.t. Fix($R_{\tilde{\mathcal{B}}}R_{\tilde{\mathcal{A}}}$) on $\mathcal{K}$.*

*Proof.* See Appendix C. □

**Theorem 3.** *Suppose that Assumptions 1-3 and 5 hold, the design matrix $\Phi$ satisfies the inequalities in Lemma 1, and $\rho_\mu \geq \frac{2}{\sigma_1}(\frac{(\theta_1+\theta_2)^2}{4\eta} + \theta_2)$ as suggested in Lemma 4. Then the sequences $\{\boldsymbol{y}_{(k)}\}_{k\in\mathbb{N}}$ and $\{\boldsymbol{\lambda}_{(k)}\}_{k\in\mathbb{N}}$ generated by Algorithm 2 satisfy $\lim_{k\to\infty}\boldsymbol{y}_k = (\mathbf{1}_N \otimes y^*)$ and $\lim_{k\to\infty}\boldsymbol{\lambda}_k = (\mathbf{1}_N \otimes \lambda^*)$, where $y^*$ is a variational GNE to problem (1) and $(y^*, \lambda^*)$ together is a solution to the KKT conditions (5).*

*Proof.* See Appendix D. □

Altogether, Theorems 2 and 3 establish the convergence of the proposed algorithm to a v-GNE regardless of the initialization. Nevertheless, the convergence speed will depend on the initial point, the choices of $\boldsymbol{\tau}_1$ to $\boldsymbol{\tau}_4$ and $\rho_\mu$, and the topology of $\mathcal{G}$. Characterizing the convergence speed is one of our future directions.

## V. Case Study and Numerical Simulations

In this study, we evaluate the proposed algorithm's performance with a Nash-Cournot game over a network. In the Nash-Cournot game, $N$ manufacturers/players indexed by $\mathcal{N}$ are involved in producing a homogeneous commodity and competing for $m$ different markets. The maximal capacities of these $m$ markets are denoted by the vector $b \in \mathbb{R}^m_{++}$.

For each manufacturer $i$ in this network, it supplies $n_i$ markets with $x_i \in \mathbb{R}^{n_i}$ units of commodities respectively. It is subject to the global market capacity constraints $A_i x_i \leq b -$ $\sum_{j\in\mathcal{N}_{-i}} A_j x_j$. The binary full-column-rank matrix $A_i \in \mathbb{R}^{m\times n_i}$ maps each entry of the decision vector $x_i$ to one among the $m$ markets. Specifically $[A_i]_{j,k} = 1$ implies that the firm $i$ supplies market $j$ with $[x_i]_k$ units of commodities. Let $n := \sum_{i=1}^N n_i$, $x := [x_1; \cdots; x_N] \in \mathbb{R}^n$, and $A := [A_1, A_2, \ldots, A_N] \in \mathbb{R}^{m\times n}$. Then, $Ax \in \mathbb{R}^m$ denotes the total quantities of commodities delivered to these $m$ markets.

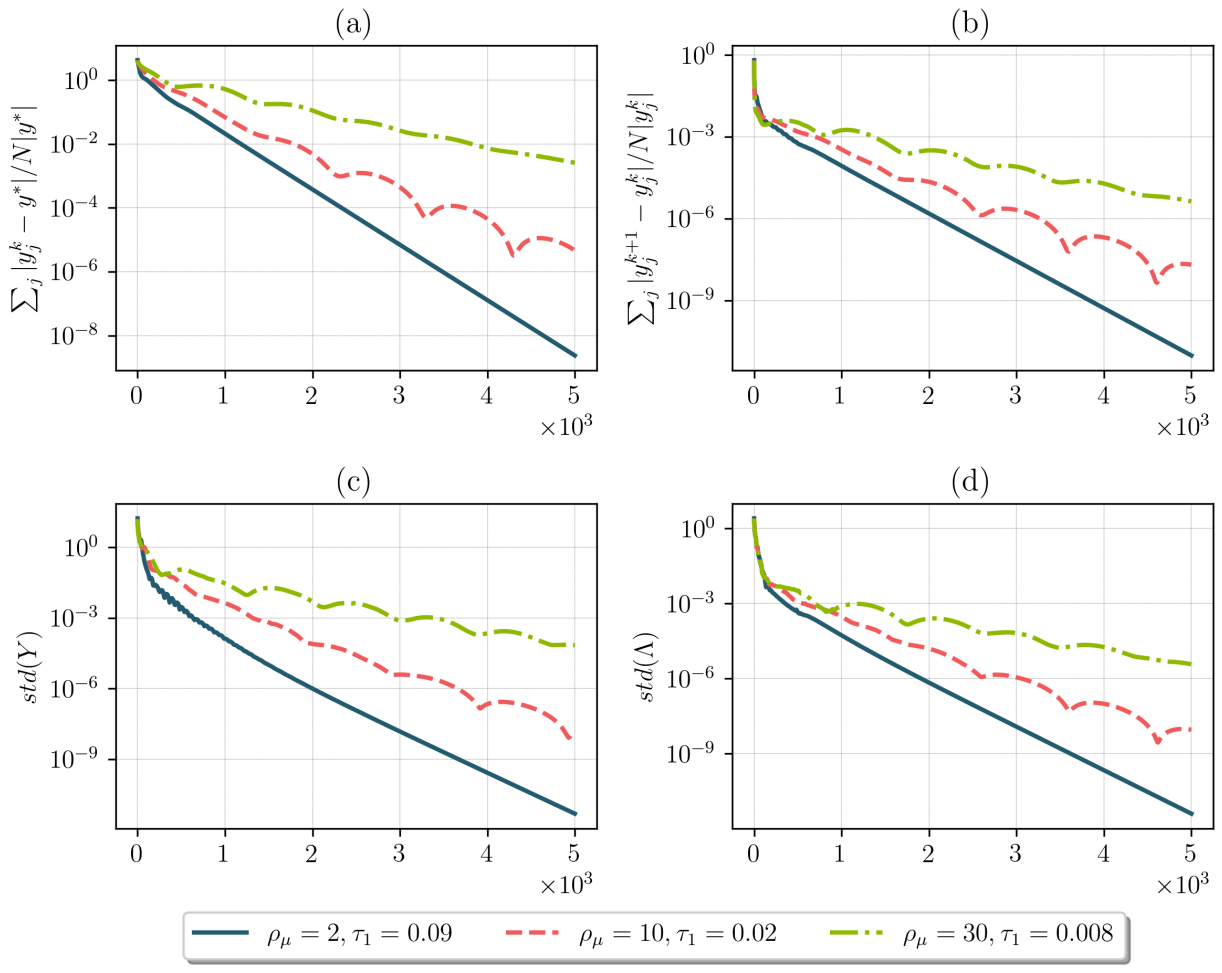

Figure 1: Convergence Result of Algorithm 2 on Nash Cournot Game (Under Assumption 4)

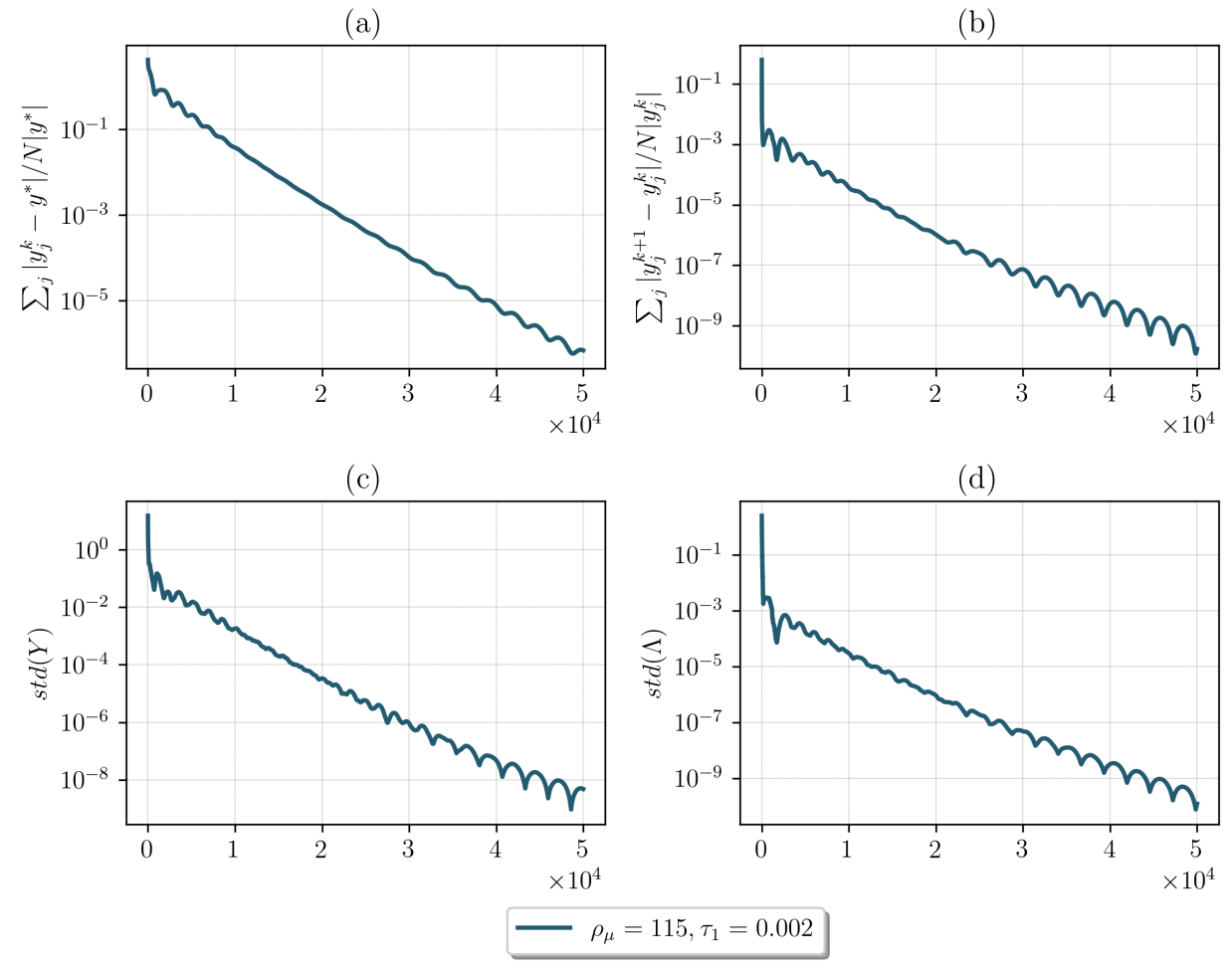

Figure 2: Convergence Result of Algorithm 2 on Nash Cournot Game (Under Assumption 5)

The local objective function is assumed to be of the form $J_i(x_i; x_{-i}) = c_i(x_i) - (P(Ax))^T A_i x_i$. We set $c_i(x_i) = x_i^T Q_i x_i + q_i^T x_i$ as the local production cost function, with $Q_i \in \mathbb{S}^{n_i}_{++}$. Let $P(u) = w - \Sigma u$ map the total quantities of supply to their unit prices, where $w \in \mathbb{R}^m_{++}$ and $\Sigma \in \operatorname{diag}(\mathbb{R}^m_{++})$.

To sum up, manufacturer $i \in \mathcal{N}$, given the supply strategies of others $(x_{-i})$, aims to solve the optimization problem (22):

$$\begin{aligned} &\underset{x_i\in\Omega_i}{\text{minimize}}\ x_i^T Q_i x_i + q_i^T x_i - (w - \Sigma\cdot Ax)^T A_i x_i \\ &\text{subject to } A_i x_i \le b - \textstyle\sum_{j\in\mathcal{N}_{-i}} A_j x_j. \end{aligned} \tag{22}$$

## A. Analysis of Nash-Cournot Game

For any $i \in \mathcal{N}$, the objective function $J_i(x_i; x_{-i})$ is a smooth function, and it satisfies the Assumption 1. To find the pseudogradient $F$, notice that $\partial_i J_i(x_i; x_{-i}) = (2Q_i + A_i^T \Sigma A_i)x_i + A_i^T \Sigma A x + q_i - A_i^T w$. By concatenating these partial derivatives together, we obtain $F: x \mapsto (M_F + A^T\Sigma A)x + \boldsymbol{q} - A^T w$, where $M_F := \operatorname{blkd}([2Q_i + A_i^T \Sigma A_i]_{i\in\mathcal{N}})$ and $\boldsymbol{q} := [q_1; \cdots; q_N]$.

Since $Q_i \succ \mathbf{0}, \forall i \in \mathcal{N}$, and $\Sigma \succ \mathbf{0}$, we can conclude that $M_F \succ 0$ and $A^T\Sigma A \succcurlyeq \mathbf{0}$. Altogether, $M_F + A^T\Sigma A \succ \mathbf{0}$. Its minimal eigenvalue $\sigma_{F\min} > 0$ and maximal eigenvalue $\sigma_{F\max} > 0$ are the strongly monotone constant ($\eta$) and the Lipschitz constant ($\theta_1$) of the pseudogradient $F$, respectively.

The extended pseudogradient $\mathcal{F}$ can be expressed as: $\mathcal{F}: \boldsymbol{x} \mapsto M_F \cdot \mathcal{R}\boldsymbol{x} + \operatorname{blkd}([A_i^T\Sigma A]_{i\in\mathcal{N}})\boldsymbol{x} + \boldsymbol{q} - A^T w$, where $\mathcal{R}$ is defined in (14). It can be further simplified to the form: $\mathcal{F}: \boldsymbol{x} \mapsto \mathcal{R}(I_N \otimes (M_F + A^T\Sigma A))\boldsymbol{x} + \boldsymbol{q} - A^T w$. The Lipschitz constant of $\mathcal{F}$ is given by the greatest singular value of $\mathcal{R}(I_N \otimes (M_F + A^T\Sigma A))$, and denoted by $\sigma_{\mathcal{F}\max}$.

## B. Simulation Results

In the numerical study, let $N = 20$ and $m = 10$. The communication graph consists of a directed circle and 10 randomly selected edges which satisfies Assumption 3. The related parameters are drawn uniformly randomly from suitable intervals. Each entry of vector $b$ satisfies $b_i \sim U[0.5, 1]$; for $P(Ax)$, each entry of vector $w$ satisfies $w_i \sim U[2, 4]$ and each diagonal entry of $\Sigma$ satisfies $\Sigma_{ii} \sim U[0.5, 0.7]$; for $c_i(x_i)$, assuming $Q_i$ is diagonal, the diagonal entry $[Q_i]_{jj} \sim U[1, 1.5]$, and each entry of $q_i$ has $[q_i]_j \sim U[0.1, 0.6]$; the feasible set $\Omega_i$ is the direct product of $n_i$ connected compact interval $[0, \Omega_{ij\max}]$, and $\Omega_{ij\max} \sim [0.2, 0.5]$; manufacturer $i$ supplies to $n_i$ markets with $n_i \sim \{2, \ldots, 6\}$.

For the case of Assumption 4, since $\mathcal{F}$ is affine in this case, it suffices to select a $\rho_\mu \ge 2$ such that $\frac{1}{2}(M_{\mathcal{F}} + M_{\mathcal{F}}^T) + \frac{\rho_\mu}{2}(L \otimes I_n) \succcurlyeq 0$, where $M_{\mathcal{F}} = \mathcal{R}^T \cdot \mathcal{R}(I_N \otimes (M_F + A^T\Sigma A))$; hence $\mathcal{R}^T\mathcal{F} + \frac{\rho_\mu}{2} L_n$ is maximally monotone. Let $\boldsymbol{\tau}_1 = \tau_1 \otimes I_{Nn}$, $\boldsymbol{\tau}_2 = \tau_2 \otimes I_{Nm}$, $\boldsymbol{\tau}_3 = \tau_3 \otimes I_{En}$, and $\boldsymbol{\tau}_4 = \tau_4 \otimes I_{Em}$, which are chosen according to Lemma 1.

For the case of Assumption 5, for $F$, $\eta = \sigma_{F\min} \approx 2.6513$, and $\theta_1 = \sigma_{F\max} \approx 10.6646$; for $\mathcal{F}$, $\theta_2 = \sigma_{\mathcal{F}\max} \approx 4.7084$. The nodes on $\mathcal{G}$ have maximal degree equal to 4. Moreover, for the Laplacian matrix $L$, $\sigma_1 \approx 0.4701$. Then, based on Lemmas 1 and 4, select $\rho_\mu = 115 \ge \frac{2}{\sigma_1}(\frac{(\theta_1+\theta_2)^2}{4\eta} + \theta_2) \approx 114.8432$, $\rho_z = 1$, $\tau_1 = 0.002$, $\tau_2 = 0.1$, and $\tau_3 = \tau_4 = 0.5$.

The performances of our proposed algorithm are illustrated in Fig.1 and 2, under the Assumptions 4 and 5, respectively. Fig.1/2a show the average of the normalized distances to the v-GNE calculated using the centralized method from [19]. Note that $y_j^k$ denotes player $j$'s local estimate of the decision vector at the $k$th iteration, and $y^*$ the generalized Nash equilibrium of the game. Fig.1/2b show the relative length of the updating step at each iteration. Let $\bar{y}^k := \frac{1}{N}\sum_{j\in\mathcal{N}} y_j^k$. Fig.1/2c exhibit how the sum of the standard deviations of the local estimates $\{y_j\}_{j\in\mathcal{N}}$, i.e., $\sum_{\ell=1}^n (\frac{1}{N}\sum_{j\in\mathcal{N}}([y_j^k]_\ell - [\bar{y}^k]_\ell)^2)^{\frac{1}{2}}$, evolves over the iterations. It measures the level of consensus among different local estimates $y_j$. Fig.1/2d are almost the same as Fig.1/2c except that we are now investigating the consensus of local dual variables $\{\lambda_j\}_{j\in\mathcal{N}}$. The numerical results verify Theorems 2 and 3 and show a linear convergence rate for each metric considered.

## VI. Conclusion and Future Directions

This paper focuses on the GNEP with generic interdependence inside the local objectives and affine coupling constraints. A distributed algorithm is proposed, which ensures exact convergence to a v-GNE and only requires local communications. For the future directions, it would be interesting to develop a set of equivalent transformations that can simplify the problem solution while preserving the v-GNE. Moreover, in this paper, we let each player keep a local copy of all players' decisions. It would be beneficial to explore the possibility of reducing the amount of local estimates for those games with structured interdependency, e.g., average/network aggregate games.

## Appendix

### A. Proof of Theorem 1

*Proof.* Suppose there exists $\omega^* \in \mathrm{Zer}(\mathcal{T})$. We then have $(B^T \otimes I_n)\boldsymbol{y}^* = \mathbf{0}$ and $(B^T \otimes I_m)\boldsymbol{\lambda}^* = \mathbf{0}$, which imply $\boldsymbol{y}^* = (\mathbf{1}_N \otimes y^*)$ and $\boldsymbol{\lambda}^* = (\mathbf{1}_N \otimes \lambda^*)$, and thus $\rho_\mu(L \otimes I_n)\boldsymbol{y}^* = \mathbf{0}$ and $(L \otimes I_m)\boldsymbol{\lambda}^* = \mathbf{0}$. Furthermore,

$$\begin{aligned}&\mathbf{0}_{mN} \in N_{\mathbb{R}^{mN}_+}(\boldsymbol{\lambda}^*) - \Lambda\mathcal{R}\boldsymbol{y}^* + \boldsymbol{b} + (B \otimes I_m)\boldsymbol{z}^*\\ &\Rightarrow \mathbf{0}_m \in (\mathbf{1}_N^T \otimes I_m)(N_{\mathbb{R}^{mN}_+}(\boldsymbol{\lambda}^*) - \Lambda\mathcal{R}\boldsymbol{y}^* + \boldsymbol{b} + (B \otimes I_m)\boldsymbol{z}^*)\\ &\Leftrightarrow \mathbf{0}_m \in N_{\mathbb{R}^m_+}(\lambda^*) - Ay^* + b \Leftrightarrow \lambda^* \geqslant \mathbf{0}_m, Ay^* - b \leqslant \mathbf{0}_m.\end{aligned}$$

Consider the first entry of the operator $\mathcal{T}$,

$$\begin{aligned}&\mathbf{0}_{nN} \in \mathcal{R}^T(\mathcal{F}(\boldsymbol{y}^*) + \Lambda^T\boldsymbol{\lambda}^*) + (B \otimes I_n)\boldsymbol{\mu}^* + N_{\tilde{\Omega}}(\boldsymbol{y}^*)\\ &\Rightarrow \mathbf{0}_n \in (\mathbf{1}_N^T \otimes I_n)(\mathcal{R}^T(\mathcal{F}(\boldsymbol{y}^*) + \Lambda^T\boldsymbol{\lambda}^*) + (B \otimes I_n)\boldsymbol{\mu}^* + N_{\tilde{\Omega}}(\boldsymbol{y}^*))\\ &\Leftrightarrow \mathbf{0}_n \in F(y^*) + A^T\lambda^* + N_\Omega(y^*)\end{aligned}$$

Hence, $(y^*, \lambda^*)$ satisfies the KKT conditions of a v-GNE (5).

In the other direction, assume the problem (5) has a v-GNE $(y^\dagger, \lambda^\dagger)$. Then for $\boldsymbol{y}^\dagger = (\mathbf{1}_N \otimes y^\dagger)$ and $\boldsymbol{\lambda}^\dagger = (\mathbf{1}_N \otimes \lambda^\dagger)$, apparently $(B^T \otimes I_n)\boldsymbol{y}^\dagger = 0$ and $(B^T \otimes I_m)\boldsymbol{\lambda}^\dagger = 0$. To prove there exists $\boldsymbol{z}^\dagger \in \mathbb{R}^{mE}$ such that $\partial_{\boldsymbol{\lambda}}\mathcal{L}^*_{\text{net}}(\boldsymbol{y}^\dagger, \boldsymbol{\lambda}^\dagger, \boldsymbol{\mu}^\dagger, \boldsymbol{z}^\dagger) \ni \mathbf{0}_{mN}$, use the second KKT condition of a v-GNE:

$$\mathbf{0}_m \le \lambda^\dagger \perp b - Ay^\dagger \ge \mathbf{0}_m. \tag{23}$$

Let $c := b - Ay^\dagger \ge \mathbf{0}_m$ and $\mathbf{0}_m \le \lambda^\dagger \perp c \ge \mathbf{0}_m$. Notice that $(B \otimes I_m)z$ spans the null space of $(\mathbf{1}_N^T \otimes I_n)$. Hence, we can find a $\boldsymbol{z}^\dagger$, such that $[b_i - A_i y_i^\dagger - \frac{1}{N}c]_{i\in\mathcal{N}} + (B \otimes I_m)\boldsymbol{z}^\dagger = \mathbf{0}_{Nm}$. It implies that $[b_i - A_i y_i^\dagger]_{i\in\mathcal{N}} + (B \otimes I_m)\boldsymbol{z}^\dagger = \frac{1}{N}(\mathbf{1}_N \otimes c)$. Based on (23), it immediately follows that

$$\mathbf{0}_{mN} \le \boldsymbol{\lambda}^\dagger \perp \boldsymbol{b} - \Lambda\mathcal{R}\boldsymbol{y}^\dagger + (B \otimes I_m)\boldsymbol{z}^\dagger \ge \mathbf{0}_{mN}. \tag{24}$$

To prove $\exists \boldsymbol{\mu}^\dagger \in \mathbb{R}^{nE}$ such that $\partial_{\boldsymbol{y}}\mathcal{L}^*_{\text{net}}(\boldsymbol{y}^\dagger, \boldsymbol{\lambda}^\dagger, \boldsymbol{\mu}^\dagger, \boldsymbol{z}^\dagger) \ni \mathbf{0}_{nN}$, observe that

$$\begin{aligned}&\mathbf{0}_{nN} \in \mathcal{R}^T(\mathcal{F}(\boldsymbol{y}^\dagger) + \Lambda^T\boldsymbol{\lambda}^\dagger) + (B \otimes I_n)\boldsymbol{\mu}^\dagger + N_{\tilde{\Omega}}(\boldsymbol{y}^\dagger) \Leftrightarrow\\ &\mathbf{0}_{nN} \in \mathcal{R}^T\big([\partial_{x_i}J_i(y^\dagger)]_{i\in\mathcal{N}} + [A_i^T\lambda^\dagger]_{i\in\mathcal{N}} + N_\Omega(y^\dagger)\big)\\ &\qquad + (B \otimes I_n)\boldsymbol{\mu}^\dagger.\end{aligned}$$

The first KKT condition of a v-GNE shows that $F(y^\dagger) + A^T\lambda^\dagger + N_\Omega(y^\dagger) \ni \mathbf{0}_n$, and thus $(\mathbf{1}_N^T \otimes I_n) \cdot \mathcal{R}^T(F(y^\dagger) + A^T\lambda^\dagger + N_\Omega(y^\dagger)) \ni \mathbf{0}_n$. Given $(B \otimes I_n)\boldsymbol{\mu}^\dagger$ spans the space $\mathrm{Null}(\mathbf{1}_N^T \otimes I_n)$, it is clear that such a $\boldsymbol{\mu}^\dagger$ exists. Altogether, $(\boldsymbol{y}^\dagger, \boldsymbol{\lambda}^\dagger, \boldsymbol{\mu}^\dagger, \boldsymbol{z}^\dagger)$ is a zero of $\mathcal{T}$. □

### B. Proof of Lemma 4

*Proof.* $\forall \omega = [\boldsymbol{y}; \boldsymbol{\lambda}; \boldsymbol{\mu}; \boldsymbol{z}]$, and $\omega^* = [\boldsymbol{y}^*; \boldsymbol{\lambda}^*; \boldsymbol{\mu}^*; \boldsymbol{z}^*] \in \mathrm{Zer}(\mathcal{T})$, $\langle \omega - \omega^*, \mathcal{A}_y(\omega) - \mathcal{A}_y(\omega^*)\rangle = \langle \boldsymbol{y} - \boldsymbol{y}^*, \mathcal{R}^T(\mathcal{F}(\boldsymbol{y}) - \mathcal{F}(\boldsymbol{y}^*))\rangle + \langle \boldsymbol{y} - \boldsymbol{y}^*, \frac{\rho_\mu}{2}(L \otimes I_n)(\boldsymbol{y} - \boldsymbol{y}^*)\rangle + \langle \boldsymbol{\lambda} - \boldsymbol{\lambda}^*, \frac{\rho_z}{2}(L \otimes I_m)(\boldsymbol{\lambda} - \boldsymbol{\lambda}^*)\rangle$. Decompose $\boldsymbol{y}$ into two orthogonal components, i.e., $\boldsymbol{y} = \boldsymbol{y}^\parallel + \boldsymbol{y}^\perp$, where $\boldsymbol{y}^\parallel = \frac{1}{N}(\mathbf{1}_{N\times N} \otimes I_n)\boldsymbol{y}$ and $\boldsymbol{y}^\perp = (I_{nN} - \frac{1}{N}(\mathbf{1}_{N\times N} \otimes I_n))\boldsymbol{y}$. Since $\omega^* \in \mathrm{Zer}(\mathcal{T})$, $\boldsymbol{y}^* = \boldsymbol{y}^{*\parallel}$. Let $\{\sigma_0, \sigma_1, \cdots, \sigma_N\}$ be the set of eigenvalues of the Laplacian matrix $L$ sorted in increasing order. By Assumption 3, we have $\sigma_0 = 0$ and $\sigma_1 > 0$. As a result, $\langle \boldsymbol{\lambda} - \boldsymbol{\lambda}^*, \frac{\rho_z}{2}(L \otimes I_m)(\boldsymbol{\lambda} - \boldsymbol{\lambda}^*)\rangle \ge 0$, $\langle \boldsymbol{y} - \boldsymbol{y}^*, \frac{\rho_\mu}{2}(L \otimes I_n)(\boldsymbol{y} - \boldsymbol{y}^*)\rangle \ge \frac{\sigma_1\rho_\mu}{2}\|\boldsymbol{y}^\perp\|^2$, and

$$\begin{aligned}&\langle \boldsymbol{y} - \boldsymbol{y}^*, \mathcal{R}^T(\mathcal{F}(\boldsymbol{y}) - \mathcal{F}(\boldsymbol{y}^*))\rangle\\ &= \langle \boldsymbol{y}^\perp + \boldsymbol{y}^\parallel - \boldsymbol{y}^*, \mathcal{R}^T(\mathcal{F}(\boldsymbol{y}) - \mathcal{F}(\boldsymbol{y}^\parallel) + \mathcal{F}(\boldsymbol{y}^\parallel) - \mathcal{F}(\boldsymbol{y}^*))\rangle\\ &= \langle \boldsymbol{y}^\perp, \mathcal{R}^T(\mathcal{F}(\boldsymbol{y}) - \mathcal{F}(\boldsymbol{y}^\parallel))\rangle + \langle \boldsymbol{y}^\parallel - \boldsymbol{y}^*, \mathcal{R}^T(\mathcal{F}(\boldsymbol{y}) - \mathcal{F}(\boldsymbol{y}^\parallel))\rangle\\ &\quad + \langle \boldsymbol{y}^\perp, \mathcal{R}^T(\mathcal{F}(\boldsymbol{y}^\parallel) - \mathcal{F}(\boldsymbol{y}^*))\rangle + \langle \boldsymbol{y}^\parallel - \boldsymbol{y}^*, \mathcal{R}^T(\mathcal{F}(\boldsymbol{y}^\parallel) - \mathcal{F}(\boldsymbol{y}^*))\rangle\\ &\ge \frac{\eta}{N}\|\boldsymbol{y}^\parallel - \boldsymbol{y}^*\|^2 - \frac{\theta_1 + \theta_2}{\sqrt{N}}\|\boldsymbol{y}^\perp\|\|\boldsymbol{y}^\parallel - \boldsymbol{y}^*\| - \theta_2\|\boldsymbol{y}^\perp\|^2.\end{aligned}$$

Merging the above three inequalities, we obtain:

$$\begin{aligned}&\langle \omega - \omega^*, \mathcal{A}_y(\omega) - \mathcal{A}_y(\omega^*)\rangle\\ &\ge \begin{bmatrix}\|\boldsymbol{y}^\parallel - \boldsymbol{y}^*\|\\ \|\boldsymbol{y}^\perp\|\end{bmatrix}^T \begin{bmatrix}\frac{\eta}{N} & -\frac{\theta_1+\theta_2}{2\sqrt{N}}\\ -\frac{\theta_1+\theta_2}{2\sqrt{N}} & \frac{\sigma_1\rho_\mu}{2} - \theta_2\end{bmatrix} \begin{bmatrix}\|\boldsymbol{y}^\parallel - \boldsymbol{y}^*\|\\ \|\boldsymbol{y}^\perp\|\end{bmatrix}.\end{aligned} \tag{25}$$

Hence, $\langle \omega - \omega^*, \mathcal{A}_y(\omega) - \mathcal{A}_y(\omega^*)\rangle \ge 0$ always holds when $\rho_\mu \ge \frac{2}{\sigma_1}(\frac{(\theta_1+\theta_2)^2}{4\eta} + \theta_2)$. □

### C. Proof of Lemma 5

*Proof.* The normal cones $N_{\tilde{\Omega}}(\boldsymbol{y})$ and $N_{\mathbb{R}^{mN}_+}(\boldsymbol{\lambda})$ are maximally monotone, since $\tilde{\Omega}$ and $\mathbb{R}^{mN}_+$ are nonempty, closed, and convex [17, Prop. 20.23, Ep. 20.26]. The skew-symmetric matrix $D$ is maximally monotone. The domain of these two operators has nonempty interior, and thus their sum is still maximally monotone. From the conclusion of Lemma 4, $\mathcal{A}_y$ is restricted monotone w.r.t. $\mathrm{Zer}(\mathcal{T})$. Since $\mathcal{A}$ is the sum of a restricted monotone operator and a monotone operator, it is restricted monotone w.r.t. $\mathrm{Zer}(\mathcal{T})$. For every $(\omega, \psi) \in \mathrm{gra}\,\tilde{\mathcal{A}}$ and every $(\omega^*, \psi^*) \in \mathrm{gra}\,\tilde{\mathcal{A}}$ with $\omega^* \in \mathrm{Zer}(\mathcal{T})$, $\Phi\psi \in \Phi\tilde{\mathcal{A}}\omega = \mathcal{A}\omega$ and $\Phi\psi^* \in \Phi\tilde{\mathcal{A}}\omega^* = \mathcal{A}\omega^*$, and therefore by restricted monotonicity of $\mathcal{A}$,

$$\langle \omega - \omega^*, \psi - \psi^*\rangle_{\mathcal{K}} = \langle \omega - \omega^*, \Phi\psi - \Phi\psi^*\rangle \ge 0. \tag{26}$$

Thus, $\tilde{\mathcal{A}}$ is restricted monotone w.r.t. $\mathrm{Zer}(\mathcal{T})$ on $\mathcal{K}$.

For every $\omega$ and $\omega^\dagger \in \mathrm{Fix}(R_{\tilde{\mathcal{B}}}R_{\tilde{\mathcal{A}}})$, we have $\omega - J_{\tilde{\mathcal{A}}}(\omega) \in \tilde{\mathcal{A}}(J_{\tilde{\mathcal{A}}}(\omega))$ and $\omega^\dagger - J_{\tilde{\mathcal{A}}}(\omega^\dagger) \in \tilde{\mathcal{A}}(J_{\tilde{\mathcal{A}}}(\omega^\dagger))$. From [17, Prop. 26.1(iii)], $\omega^\dagger \in \mathrm{Fix}(R_{\tilde{\mathcal{B}}}R_{\tilde{\mathcal{A}}})$ implies $J_{\tilde{\mathcal{A}}}(\omega^\dagger) \in \mathrm{Zer}(\mathcal{T})$. By the restricted monotonicity of $\tilde{\mathcal{A}}$ shown above,

$$\langle (\omega - J_{\tilde{\mathcal{A}}}(\omega)) - (\omega^\dagger - J_{\tilde{\mathcal{A}}}(\omega^\dagger)), J_{\tilde{\mathcal{A}}}(\omega) - J_{\tilde{\mathcal{A}}}(\omega^\dagger)\rangle_{\mathcal{K}} \ge 0, \tag{27}$$

which proves the restricted firm nonexpansiveness of $J_{\tilde{\mathcal{A}}}$ w.r.t. $\mathrm{Fix}(R_{\tilde{\mathcal{B}}}R_{\tilde{\mathcal{A}}})$ on $\mathcal{K}$ (see [17, Prop. 4.4]). The restricted nonexpansiveness of $R_{\tilde{\mathcal{A}}}$ directly follows from this result. □

### D. Proof of Theorem 3

*Proof.* The DR algorithm follows the Krasnosel'skii-Mann iteration scheme $\tilde{\omega}^{(k+1)} = \tilde{\omega}^{(k)} + \frac{1}{2}(R_{\tilde{\mathcal{B}}}R_{\tilde{\mathcal{A}}} - I)\tilde{\omega}^{(k)}$, and generates a well-defined sequence. For notational brevity, let $\tilde{R} := R_{\tilde{\mathcal{B}}}R_{\tilde{\mathcal{A}}}$. Then we have, $\forall \tilde{\omega}^* \in \mathrm{Fix}(\tilde{R})$,

$$\begin{aligned}&\left\|\tilde{\omega}^{(k+1)} - \tilde{\omega}^*\right\|_{\mathcal{K}}^2 = \left\|\frac{1}{2}(\tilde{\omega}^{(k)} - \tilde{\omega}^*) + \frac{1}{2}(\tilde{R}\tilde{\omega}^{(k)} - \tilde{\omega}^*)\right\|_{\mathcal{K}}^2 \\ &= \frac{1}{2}\left\|\tilde{\omega}^{(k)} - \tilde{\omega}^*\right\|_{\mathcal{K}}^2 + \frac{1}{2}\left\|\tilde{R}\tilde{\omega}^{(k)} - \tilde{\omega}^*\right\|_{\mathcal{K}}^2 - \frac{1}{4}\left\|\tilde{R}\tilde{\omega}^{(k)} - \tilde{\omega}^{(k)}\right\|_{\mathcal{K}}^2 \\ &\le \left\|\tilde{\omega}^{(k)} - \tilde{\omega}^*\right\|_{\mathcal{K}}^2 - \frac{1}{4}\left\|\tilde{R}\tilde{\omega}^{(k)} - \tilde{\omega}^{(k)}\right\|_{\mathcal{K}}^2,\end{aligned} \tag{28}$$

where the inequality follows from the nonexpansiveness of $R_{\tilde{\mathcal{B}}}$ as well as the restricted nonexpansiveness of $R_{\tilde{\mathcal{A}}}$. Therefore, $\{\tilde{\omega}^{(k)}\}_{k\in\mathbb{N}}$ is Fejer monotone w.r.t. $\mathrm{Fix}(R_{\tilde{\mathcal{B}}}R_{\tilde{\mathcal{A}}})$.

We can derive from (28) that $\sum_{k\in\mathbb{N}} \frac{1}{4}\left\|\tilde{R}\tilde{\omega}^{(k)} - \tilde{\omega}^{(k)}\right\|_{\mathcal{K}}^2 \le \left\|\tilde{\omega}^{(0)} - \tilde{\omega}^*\right\|_{\mathcal{K}}^2$. Therefore, $\lim_{k\to\infty}\left\|\tilde{R}\tilde{\omega}^{(k)} - \tilde{\omega}^{(k)}\right\|_{\mathcal{K}}^2 = 0$.

From the Fejer monotonicity of $\{\tilde{\omega}^{(k)}\}_{k\in\mathbb{N}}$, it is a bounded sequence and thus has a convergent subsequence $\{\tilde{\omega}^{(k_i)}\}_{i\in\mathbb{N}}$ such that $\lim_{i\to\infty}\tilde{\omega}^{(k_i)} = \tilde{\omega}^\dagger$. Since $\tilde{R}$ is a continuous mapping (see Appendix E), $\lim_{i\to\infty}\tilde{R}\tilde{\omega}^{(k_i)} - \tilde{\omega}^{(k_i)} = 0$ implies $\tilde{R}\tilde{\omega}^\dagger = \tilde{\omega}^\dagger$, and hence $\tilde{\omega}^\dagger \in \mathrm{Fix}(\tilde{R})$. This implies that the distance between $\{\tilde{\omega}^{(k)}\}_{k\in\mathbb{N}}$ and $\tilde{\omega}^\dagger$ satisfies

$$\liminf_{k\to\infty} d_{\tilde{\omega}^\dagger}(\tilde{\omega}^{(k)}) = 0. \tag{29}$$

Since $\tilde{\omega}^\dagger \in \mathrm{Fix}(R_{\tilde{\mathcal{B}}}R_{\tilde{\mathcal{A}}})$, we have $\left\|\tilde{\omega}^{(k+1)} - \tilde{\omega}^\dagger\right\|_{\mathcal{K}}^2 \le \left\|\tilde{\omega}^{(k)} - \tilde{\omega}^\dagger\right\|_{\mathcal{K}}^2 - \frac{1}{4}\left\|R_{\tilde{\mathcal{B}}}R_{\tilde{\mathcal{A}}}\tilde{\omega}^{(k)} - \tilde{\omega}^{(k)}\right\|_{\mathcal{K}}^2$ which implies that $\{d_{\tilde{\omega}^\dagger}(\tilde{\omega}^{(k)})\}_{k\in\mathbb{N}}$ is a monotonically non-increasing sequence. Combining (29), we have $d_{\tilde{\omega}^\dagger}(\tilde{\omega}^{(k)}) \to 0$, and hence $\{\tilde{\omega}^{(k)}\}_{k\in\mathbb{N}}$ converges to the point $\tilde{\omega}^\dagger$. Altogether, $\omega^\dagger = J_{\tilde{\mathcal{A}}}\tilde{\omega}^\dagger$ is a zero of $\mathcal{T}$. Using the results from Theorem 1, the proof is complete. □

### E. The Continuity of the Composite Operator

**Lemma 6.** *Suppose Assumptions 1, 2 and 5 hold. Then the composite operator $R_{\tilde{\mathcal{B}}}R_{\tilde{\mathcal{A}}}$ is continuous.*

*Proof.* By Lemma 2, $R_{\tilde{\mathcal{B}}}$ is nonexpansive on $\mathcal{K}$, and hence Lipschitz continuous on $\mathcal{K}$. Since the composition of two continuous operators is still continuous and $R_{\tilde{\mathcal{A}}} = 2J_{\tilde{\mathcal{A}}} - I$, it is sufficient to prove that $J_{\tilde{\mathcal{A}}}$ is continuous. The explicit formulas of $J_{\tilde{\mathcal{A}}}$ inside Algorithm 2 imply that $\boldsymbol{\mu}$ and $\boldsymbol{z}$ are updated using simple affine functions. Moreover, the update of $\boldsymbol{\lambda}$ is an affine function followed by a projection onto the positive orthant, and thus is also continuous. We can split the resolvent w.r.t. $\boldsymbol{y}$ into two decoupled argmin problems w.r.t. $y_i^i$ and $y_i^{-i}$, $\forall i \in \mathcal{N}$. Observe that the update of $y_i^{-i}$ is also continuous w.r.t. $\tilde{\omega}$. The update of $y_i^i$ is given as follows:

$$\begin{aligned}y_i^i = \operatorname*{argmin}_{v_i\in\Omega_i}\Big[&J_i(v_i; y_i^{-i}) + \frac{1}{2}\big(\textstyle\sum_{j\in\mathcal{N}_i^+}\tilde{\mu}_{ji}^i - \sum_{j\in\mathcal{N}_i^-}\tilde{\mu}_{ij}^i\big)^T v_i \\ &+ \frac{1}{2}\tilde{\lambda}_i^T A_i v_i + \frac{\rho\mu}{2}\textstyle\sum_{j\in\mathcal{N}_i}(\tilde{y}_i^i - \tilde{y}_j^i)^T v_i + \frac{1}{2\tau_{1i}}\left\|v_i - \tilde{y}_i^i\right\|^2\Big].\end{aligned} \tag{30}$$

Denote the part inside the square brackets of (30) by $\hat{J}_i(v_i; y_i^{-i}, \tilde{\omega})$. Based on Assumption 5, $\partial_{v_i}\hat{J}_i$ is a continuous mapping in terms of $v_i$, $y_i^{-i}$, and $\tilde{\omega}$. Given any fixed $y_i^{-i}$ and $\tilde{\omega}$, $\hat{J}_i(v_i; y_i^{-i}, \tilde{\omega})$ is a strongly convex function. Thus $\partial_{v_i}\hat{J}_i(\cdot; y_i^{-i}, \tilde{\omega}) : \Omega_i \to \mathbb{R}$ is a one-to-one mapping. Applying the generalization of Implicit Function Theorem [20, Th. 1.1], the solution $y_i^i$ to the equation $\partial_{v_i}\hat{J}_i(v_i; y_i^{-i}, \tilde{\omega}) = 0$ is continuous w.r.t. $y_i^{-i}$ and $\tilde{\omega}$. □